\documentclass[a6paper,12pt]{article} 
\newdimen\paperhight
\usepackage{amsmath,amssymb}
\usepackage[usenames]{color}
\usepackage{graphicx}
\usepackage{amsfonts}
\usepackage{ascmac}

\newcommand{\dsp}{\displaystyle}

\newcommand{\qed}{\quad\hbox{\rule[-2pt]{3pt}{6pt}}\par\vspace{5mm}}
\newcommand{\pr}{\par \vspace{3mm} \noindent {\bf [Proof]} \qquad}
\newcommand{\prend}{\hfill \qed}

\newcommand{\1}{{\bf 1}} 
\newcommand{\C}{\mathbb C} 
\newcommand{\Z}{\mathbb Z} 
\newcommand{\Q}{\mathbb Q} 
\newcommand{\N}{\mathbb N}

\newcommand{\CF}{{\cal F}}

\newcommand{\CR}{{\cal R}}

\newcommand{\hf}{{\frac{1}{2}}}
\newcommand{\st}{{\frac{1}{16}}}

\newcommand{\wt}{{\rm wt}}
\newcommand{\Hom}{{\rm Hom}}
\newcommand{\End}{{\rm End}}

\newtheorem{thm}{Theorem}
\newtheorem{prn}[thm]{Proposition}

\newtheorem{cry}[thm]{Corollary}
\newtheorem{lmm}[thm]{Lemma}
\newtheorem{rmk}{Remark}

\begin{document}
\title{$C_2$-cofiniteness of cyclic-orbifold models}
\author{\begin{tabular}{c}
Masahiko MIYAMOTO\\
Institute of Mathematics, \\
University of Tsukuba, \\
Tsukuba, 305 Japan \end{tabular}}
\date{}
\maketitle

\begin{abstract}
We prove an orbifold conjecture for a solvable automorphism group. Namely, 
we show that if $V$ is a $C_2$-cofinite simple vertex operator algebra and 
$G$ is a finite solvable automorphism group of $V$, 
then the fixed point vertex operator subalgebra $V^G$ is also $C_2$-cofinite. 
This offers a mathematically rigorous background to orbifold theories 
with solvable automorphism groups. 
\end{abstract}

\section{Introduction}
In order to explain the moonshine phenomenon on the monster simple group 
and the modular functions, 
Bocherds \cite{B} has introduced a concept of vertex (operator) algebra as an 
algebraic version of conformal field theory. 
It is a quadruple $(V,Y,\1,\omega)$ satisfying the several axioms, where 
$V$ is a graded vector space $V=\oplus_{i=-K}^{\infty}V_i$, 
$Y(v,z)=\sum_{m\in \Z} v_mz^{-m-1}\in \End(V)[[z,z^{-1}]]$ denotes 
a vertex operator of $v\in V$ on $V$ which satisfies Borcherds identity (2.1),  
$\1\in V_0$ and $\omega\in V_2$ are specified elements called 
the vacuum and the Virasoro element of $V$, respectively. We set 
$Y(\omega,z)=\sum_{n\in \Z} L(n)z^{-n-2}$. 

One of the main targets in the research of vertex operator algebras (shortly VOA) 
is a construction of VOAs of finite type, that is, all modules (including weak modules) 
have a composition series consisting of only finitely many isomorphism classes of simple modules. 
If $V$ is a VOA and $\sigma$ is a finite automorphism of order $p$, then 
a fixed point subVOA $V^\sigma$ is called an orbifold model, (see \cite{Dix}, \cite{Dij}). 
So-called ``orbifold conjecture'' says that if $V$ is of finite type, then so is $V^\sigma$. 
It is revealed that the above finiteness condition is equivalent to 
the $C_2$-cofiniteness by \cite{Bu} and \cite{M1}. 
Here a $V$-module $W$ is called to be $C_2$-cofinite when 
$C_2(W)={\rm Span}_{\C}\{v_{-2}u \mid v\in V, u\in W, \wt(v)>0\}$ has a finite 
co-dimension in $W$. This condition was originally introduced by \cite{Z} 
as a technical condition to prove the modular invariance property. 
This condition is powerful and the most general theorems require this condition. 
For example, the author in \cite{M2} mentioned that 
if the orbifold model $V^\sigma$ is $C_2$-cofinite, then we are able to get 
all information of (weak) $V^\sigma$-modules from (twisted and ordinary) $V$-modules and 
every simple $V^\sigma$-module 
is a submodule of a (twisted or ordinary) $V$-module. 
Therefore, the $C_2$-cofiniteness on $V^\sigma$ offers a mathematically rigorous background 
to all orbifold theories.

For the orbifold conjecture, there are partial answers. For example, 
T.~Abe has proved it for a permutation automorphism of order $p=2$. 
For lattice VOAs, Yamskulna \cite{Y} has shown the case $p=2$ and the author \cite{M3} 
has shown the case $p=3$, which was used to construct a new holomorphic VOA of central 
charge $24$. 
In this paper, we will prove the all cases for any finite order $p$ with the powerful help of 
Borcherds identity (2.1) and the skew-symmetry (4.1). \\

\noindent
{\bf Main Theorem} \quad 
{\it Let $V$ be a $C_2$-cofinite simple VOA of CFT-type and 
$\sigma\in {\rm Aut}(V)$ of finite order $p$. 
Then a fixed point vertex operator subalgebra 
$V^{\sigma}$ is also $C_2$-cofinite.} \\

As corollaries, we have:

\begin{thm}\label{2G} \quad
Let $V$ be a $C_2$-cofinite simple VOA of CFT-type and 
$G\leq {\rm Aut}(V)$ finite solvable. 
Then a fixed point vertex operator subalgebra 
$V^{G}$ is also $C_2$-cofinite.
\end{thm}

\begin{cry}\label{CL} \quad 
Let $V$ be a $C_2$-cofinite VOA and a subVOA $U$ is isomorphic to a lattice
VOA. Then the commutant $E$ of $U$ is $C_2$-cofinite.
\end{cry}

\begin{cry}\label{CI} \quad 
If $V$ is $C_2$-cofinite and a subVOA $U$ is isomorphic to an 
$2$-dim. Ising model $L(\hf,0)$ of central charge $\hf$, 
then the commutant $E$ of $U$ is $C_2$-cofinite.
\end{cry}

Here the commutant $E$ of $U$ is defined 
by $\{v\in V\mid u_mv=0 \mbox{ for all } u\in U, m\geq 0\}$ and it is a subVOA. 

\begin{rmk}
In this paper, we assume that $V$ is of CFT-type. 
This is because of simplifying the proof. From our proof, it is not difficult to 
see that we have the same conclusion without the assumption of CFT-type. 
\end{rmk}

We close this introduction by acknowledging with thanks a number of 
communications with Yu-ichi Tanaka and Shigeki Akiyama. The author thanks 
Toshiyuki Abe, Hiroshi Yamauchi and Atsushi Matsuo for 
reviewing the manuscript and their suggestions about the shorter proofs.  
He also thanks to the organizers of the conference 
held at Taitung university in March 2013 for their hospitality.

\section{Truncation property} 
From the axiom of VOAs, for $v\in V_r$ and $u\in V_n$, we have 
$v_mu\in V_{r-m-1+n}$. Hence 
there is an integer $N$ such that $v_{n}u=0$ for any $n>N$. 
This is called a truncation property. 
To simplify the notation, we will say that $v$ is truncated on $u$.  

Set $V^{\ast}=\Hom(V,\C)$ and define a pairing $\langle \cdot,\cdot\rangle$ by 
$\langle v, \xi\rangle=\xi(v)$ for $\xi\in V^*$ and $v\in V$. 
For $v\in V$ and $m\in \Z$, actions $v_m$ on $V^{\ast}$ are defined by 
$$\langle w, Y^{\ast}(v,z)\xi \rangle=
\langle Y(e^{L(1)z}(-z^{-2})^{L(0)}v,z^{-1})w,\xi\rangle $$
for $w\in V, \xi\in \Hom(V,\C)$, where 
$Y^{\ast}(v,z)=\sum_{m\in \Z}v_mz^{-m-1}$ is called an adjoint operator of $v$. 
An important fact is that 
$(\oplus_{m=0}^{\infty}\Hom(V_m,\C),Y^{\ast})$ becomes a $V$-module, 
see \cite{FHL} for the proof. This module is called 
a restricted dual of $V$ which is denoted by $V'$. 
In particular, $Y^{\ast}(\cdot,z)$ satisfy the Borcherds identity:\\
$$\sum_{i=0}^{\infty}\binom{m}{i}(u_{r+i}v)_{m+n-i} \xi=
\sum_{j=0}^{\infty}(-1)^j\binom{r}{j}(u_{r+m-j}v_{n+j}\xi
-(-1)^rv_{r+n-j}u_{m+j}\xi) \eqno{(2.1)}$$
for any $m,n,r\in \Z$, $v,u\in V$, $\xi\in V'$. 
Since $V^{\ast}=\prod_n{\rm Hom}(V_n,\C)$,  
we can express $\xi\in V^{\ast}$ by $\prod_n \xi_{(n)}$ with $\xi_n\in \Hom(V_n,\C)$. 
We call $\xi\in V^{\ast}$ $L(0)$-free if $\dim \C[L(0)]\xi =\infty$, that is, 
$\xi_{(n)}\not=0$ for infinitely many $n$. 
We note that if $V$ is $C_2$-cofinite, then any (weak) module does not contain 
$L(0)$-free elements.

The weight of the terms in (2.1) for $\xi\in \Hom(V_t,\C)$ and that for 
$\xi\in \Hom(V_s,\C)$ are 
different when $t\not=s$. We also have that the both sides of (2.1) are 
well-defined for each $\xi\in \Hom(V_t,\C)$. Therefore  
the Borcherds' identity is also well-defined on $V^{\ast}$, 
as Haisheng Li has pointed out in \cite{L}.  
However, $V^{\ast}$ is not a $V$-module. The problem is a failure of truncation properties. 

\begin{lmm}\quad
If $u$ and $v$ are truncated on $\xi$, then $v_mu$ is also truncated on $\xi$ 
for any $m$.  In particular, if $V$ is generated by $\Omega\subseteq V$ as a vertex algebra and 
all elements in $\Omega$ are truncated on $\xi$, then all elements in $V$ are truncated on $\xi$. 
\end{lmm}

\pr 
We may assume $u_n\xi=v_n\xi=u_nv=0$ for $n\geq N$. 
We assert that for $s\in \N$ and $n\geq 2N+s$, we have $(u_{N-s}v)_n\xi=0$. 
Suppose false and let $s$ be a minimal counterexample. 
Substituting $r=N-s$, $n=N+s+p$, $m=N+q$ in (2.1) with $p,q\geq 0$, the left side equals  
$$\begin{array}{rl}
{\rm LH}=&{\dsp \sum_{i=0}^{\infty}\binom{N+q}{i}(u_{N-s+i}v)_{2N+q+s+p-i} \xi=
\sum_{i=0}^s\binom{N+q}{i}(u_{N-(s-i)}v)_{2N+s-i+p+q}\xi}\cr
=&{\dsp (u_{N-s}v)_{2N+s+p+q}\xi }\end{array} $$
by the minimality of $s$. On the other hand, the right side is  
$$
{\dsp {\rm RH}=\sum_{i=0}^{\infty}(-1)^i\binom{N\!-\!s}{i}\left(u_{2N-s+q-i}v_{N+s+p+i}\xi
-(-1)^{N-s}v_{2N-s+p-i}u_{N+q+i}\xi\right)=0,}$$
which contradicts the choice of $s$. 
\prend

Since ${\dsp v_nu_m\xi=u_mv_n\xi+\sum_{i=0}^{\infty}\binom{n}{i}(v_iu)_{n+m-i}\xi}$, 
Lemma 2 (see also H.~Li \cite{L}) implies:

\begin{lmm}\quad
If $v,u\in V$ are truncated on $\xi\in V^{\ast}$, then $v$ is truncated on $u_m\xi$ for any $m$. 
In particular, if all elements of $V$ are truncated on $\xi$, then 
${\rm Span}_{\C}\{u^1_{m_1}\cdots u^k_{m_k}\xi\mid u^i\in V, m_i\in \Z\}$ is a $V$-module. 
\end{lmm}

\section{General setting}
Let $(V,Y,\1,\omega)$ be a $C_2$-cofinite VOA 
and $\sigma$ an automorphism of $V$ of order $p$. 
Viewing $V$ as a $<\sigma>$-module, we decompose 
$$     V=V^{(0)}\oplus V^{(1)}\oplus \cdots \oplus V^{(p-1)} $$
where $V^{(m)}=\{ v\in V \mid v^\sigma=e^{2\pi \sqrt{-1}m/p}v \}$. 
For subsets $A,B$ of $V$ and $m\in \Z$, $A_{(m)}B$ denotes a subspace 
${\rm Span}_{\C}\{a_mb\mid a\in A, b\in B\}$. 

\begin{lmm}\quad
$(V^{(1)})_{(-2)}V^{(0)}+(V^{(0)})_{(-2)}V^{(1)}$ has a finite co-dimension in $V^{(1)}$. 
\end{lmm}

\pr
Suppose false, i.e. 
$V^{(1)}_m/((V^{(1)})_{(-2)}V^{(0)}+(V^{(0)})_{(-2)}V^{(1)})_m\not=0$ for infinitely many $m$. 
Then there is a $L(0)$-free element $\xi\in (V^{(1)})^{\ast}$ such that 
$\langle (V^{(1)})_{(-2)}V^{(0)}+(V^{(0)})_{(-2)}V^{(1)},\xi\rangle=0$. In other words, 
$$ 0=\langle v_{-2-\N}u,\xi\rangle=\langle u_{-2-\N}v,\xi\rangle $$
for any $v\in V^{(1)}$ and $u\in V^{(0)}$. Since $V^{(1)}$ is a direct summand of $V$, we may view 
$(V^{(1)})^{\ast}\subseteq V^{\ast}$. 
By taking adjoint operators, we have:
$$\begin{array}{l}
{\dsp \langle u, v_{2\wt(v)+\N}\xi \rangle 
=\langle (-1)^{\wt(v)} \sum_{s=0}^{\infty}\frac{1}{s!}(L(1)^sv)_{-2-s-\N}u, \xi\rangle=0  }\cr
{\dsp \langle v, u_{2\wt(u)+\N}\xi \rangle 
=\langle (-1)^{\wt(u)} 
\sum_{s=0}^{\infty}\frac{1}{s!}(L(1)^su)_{-2-s-\N}v, \xi\rangle=0, }\end{array}$$
which imply that $v\in V^{(1)}$ and $u\in V^{(0)}$ truncate on $\xi$. 
However, since $V$ is simple, $V^{(1)}+V^{(0)}$ generates a $C_2$-cofinite VOA 
$V$ by normal products, which contradicts that $\xi$ is $L(0)$-free.  
\prend

So, there is a finite dimensional subspace $P$ of $V^{(1)}$ 
such that $V^{(1)}=(V^{(1)})_{(-2)}V^{(0)}+(V^{(0)})_{(-2)}V^{(1)}+P$. 
We may assume that $P$ is a direct sum of homogeneous spaces. 

\begin{prn}\label{F1}\quad  $V^{(1)}=(V^{(0)})_{(-2)}V^{(1)}+\C [L(-1)]P$. 
\end{prn}

\pr
Suppose false and we choose 
$0\not=w\in V^{(1)}-((V^{(0)})_{(-2)}V^{(1)}+\C [L(-1)]P)$ with minimal weight. 
Since $V^{(1)}=(V^{(1)})_{(-2)}V^{(0)}+(V^{(0)})_{(-2)}V^{(1)}+P$, 
we may assume $w\in (V^{(1)})_{(-2)}V^{(0)}$.  
We may also assume $w=a_{-2}u$ with $a\in V^{(1)}$ and $u\in V^{(0)}$. 
Then by the skew-symmetry (4.1), we have 
$$w=-u_{-2}a-\sum_{j=1}^{\infty}\frac{(-1)^j}{j!}L(-1)^ju_{-2+j}a.$$
Since $\wt(u_{-2+j}a)<\wt(u_{-2}a)=\wt(w)$ for $j\geq 1$, we have 
$$w\in (V^{(0)})_{(-2)}V^{(1)}+\C[L(-1)](V^{(0)})_{(-2)}V^{(1)}+\C [L(-1)]P
\subseteq (V^{(0)})_{(-2)}V^{(1)}+\C [L(-1)]P$$ 
by the minimality of $\wt(w)$, which contradicts the choice of $w$. 
\prend

\section{The coefficient functions} 
Since $L(-1)C_2(V^{(1)})\subseteq C_2(V^{(1)})$, $V^{(1)}/C_2(V^{(1)})$ is a finitely generated 
$\C[L(-1)]$-module by Proposition \ref{F1}. 
Let $T$ be the inverse image in $V^{(1)}$ of the $L(-1)$-torsion 
submodule of $V^{(1)}/C_2(V^{(1)})$.  
Then there is a set of free generators $\{\alpha^i:i=1,...,t\}$ such that 
$$V^{(1)}=\left( \oplus_{i=1}^t \C[L(-1)]\alpha^i \right)\oplus T.$$ 
If $V^{(1)}$ and $V^{(p-1)}$ are $C_2$-cofinite, then so is $V^{(0)}$ by the main theorem 
in \cite{M4}.  So we may assume that $V^{(1)}$ is not $C_2$-cofinite, that is, $t\geq 1$. 

The key idea in this paper is to denote elements $a_{-n}b$ in $V^{(1)}$ 
as a linear combination 
$\sum_{i=1}^t f^i(n)(\alpha^i)_{-\wt(a)-\wt(b)-n+\wt(\alpha^i)}\1$ modulo $T$ 
for a sufficiently large $n$ and consider $f^i(n)$ as functions of $n\in \Z$. 
We note $\frac{L(-1)^m}{m!}\alpha=\alpha_{-1-m}\1$ for $m\in \N$.

We may choose $\alpha^1$ so that $\wt(\alpha^1)$ is 
the minimal weight of elements in $V^{(1)}-T$. From now on, 
we denote $\alpha^1$ by $\alpha$.   
Since it is enough to prove Main Theorem, we will concentrate only on the 
coefficients of $\alpha_{-m}\1$. In order to do it, 
we will use $\equiv$ to denote an equivalence relation modulo   
$$\widehat{T}:=T+\C[L(-1)]\alpha^{2}+\cdots+\C[L(-1)]\alpha^{t}.$$
Under this setting, we view $x\in \Z$ as a variable and 
consider functions $f(x)$ given by 
$$ a_{-x-\wt(\alpha)+\wt(a)+\wt(b)-1}b\equiv f(x)\alpha_{-x-1}\1 \pmod{\widehat{T}}.$$
From now on, for $a,b\in V$, we always use $M$ to denote $\wt(a)+\wt(b)-\wt(\alpha)$. 
We note that since $\wt(a_{-x+M-1}b)<\wt(\alpha)$ for $x\in \Z_{<0}$, 
$a_{-x+M-1}b\in \widehat{T}$ by the choice of $\alpha$.  
Namely, $f(x)=0$ for $x\in \Z_{<0}$. 
So we will consider the set ${\rm Map}(\N,\C)$ of 
all maps from $\Z$ to $\C$ satisfying $f(n)=0$ for $n\in \Z_{<0}$.  

For each $V^{(k)}$, we consider the space of coefficients $f(x)$ of $a_{-x+M-1}b$ at 
$\alpha_{-x-1}\1$ modulo $\widehat{T}$ for $a\in V^{(k)}$ and $b\in V^{(p+1-k)}$, that is, 
$$ \CF_k={\rm Span}_{\C}\left\{ f\in {\rm Map}(\N,\C)\mid 
\begin{array}{l}
{}^{\exists}a\in V^{(k)},{}^{\exists}b\in V^{(p-k+1)} \mbox{ s.t. }  \cr 
 a_{-x+M-1}b\equiv f(x)\alpha_{-x-1}\1 \pmod{\widehat{T}}
\mbox{ for }x\in \N  \end{array}\right\}.$$

\begin{lmm}\label{ID}\quad
$\CF_k$ are all $\C[x]$-invariant.
\end{lmm}

\pr 
Clearly, $\CF_k$ is a vector space. 
If $a_{-x+M-1}b\equiv f(x)\alpha_{-x-1}\1$, then 
$$ (L(-1)a)_{-x+M}b=(x-M)a_{-x+M-1}b\equiv (x-M)f(x)\alpha_{-x-1}\1 \pmod{\widehat{T}}. $$
Hence $xf(x)\in \CF_k$. 
\prend

\begin{lmm}\quad
For $f(x)\in \CF_0$, $Q_f=\{n\in \Z\mid f(n)\not=0\}$ is a finite set.  
\end{lmm}

\pr
Let $a\in V^{(0)}$ and $b\in V^{(1)}$ and set $a_{-x-1+M}b\equiv f(x)\alpha_{-x-1}\1 
\pmod{\widehat{T}}$. 
Then since $a_{-x+M-1}b\in C_2(V^{(1)})$ for $x\geq M+1$, 
$f(x)=0$ for $x\geq M+1$. We also have $f(x)=0$ for $x<0$ as we have shown. 
Since all elements in $\CF_0$ are finite linear combinations of 
such elements, we have the desired result.
\prend

For a map $f:\N\to \C$, we introduce two operators $S$ and $T$ as follows:
$$\begin{array}{ll}
{\dsp Sf(n)=(-1)^{n}\sum_{k=0}^n\binom{n}{k}(-1)^kf(n-k)}&\mbox{for }n\in \N \cr
{\dsp Tf(n)=(-1)^nf(n)}&\mbox{for }n\in \N. 
\end{array}$$

Clearly, $S^2=T^2={\rm id}$. We also have the following by induction.

\begin{lmm}\quad 
$$(ST)^kf(n)=\sum_{j=0}^n\binom{n}{j}k^jf(n-j) \mbox{ for }k=1,...$$
\end{lmm} 

The following is a key result. 

\begin{lmm}\label{ST}\quad 
For $k=1,2,...,p-1$, we have:
$$\CF_{1+p-k}=S(\CF_{k}) \quad\mbox{ and  }\quad
\CF_{p-k}=T(\CF_k).$$ 
\end{lmm}

\pr
The operator $S$ comes from the skew-symmetry. 
In fact, for $a\in V^{(k)}, b\in V^{(p-k+1)}$ and 
$a_{-x+M-1}b\equiv f(x)\alpha_{-x-1}\1$ for $x\in \N$, then 
$$\begin{array}{rl}
b_{-x+M-1}a\equiv& {\dsp (-1)^{x+M}\sum_{k=0}^{\infty}\frac{L(-1)^k}{k!}(-1)^ka_{-(x-k)+M-1}b 
\hspace{4cm}\hfill (4.1) }\cr 
\equiv& {\dsp (-1)^{x+M}
\sum_{k=0}^{x}\frac{L(-1)^k}{k!}(-1)^kf(x-k)\alpha_{-x+k-1}\1} \cr
 &\mbox{since $L(-1)^ka_{-x+k+M-1}b\in \widehat{T}$ for $k>x$,} \cr
\equiv&{\dsp (-1)^{x+M}\sum_{k=0}^{x}\binom{x}{k}(-1)^kf(x-k)\alpha_{-x-1}\1}\cr
\equiv&{\dsp (-1)^{M}Sf(x)\alpha_{-x-1}\1.}
\end{array}$$ 
Therefore, $\CF_{1+p-k}\subseteq S(\CF_k)$. Since $S^2=1$, we have 
the equality $\CF_{1+p-k}=S(\CF_k)$.

For any $m\in \Z$ and $h\in \N$, 
by substituting $n=-m-2-h$ and $r=-x+N+h$ in the Borcherds' identity (2.1), we have 
$$\begin{array}{l}
\dsp{ 0\equiv \sum_{i=0}^{\infty}\binom{m}{i}(u_{-x+N+h+i}v)_{-2-h-i} \xi}\cr
\dsp{\mbox{}\quad =\sum_{j=0}^{\infty}\binom{-x+N+h}{j}(-1)^ju_{-x+N+h+m-j}v_{-m-2-h+j}\xi }\cr
\dsp{\mbox{}\qquad -(-1)^{-x}\sum_{j=0}^{-m+\wt(u)+\wt(\xi)}\binom{-x+N+h}{j} (-1)^{j+N+h}v_{-x+N-m-2-j}u_{m+j}\xi)}
\end{array}$$
for $u\in V^{(k)}$, $v\in V^{(p-k)}$, $\xi\in V^{(1)}$, where $N=\wt(\xi)+\wt(u)+\wt(v)$. 
We note $u_{m+j}\xi=0$ for $j\geq Q=\wt(u)+\wt(\xi)-m$. Since we will treat only $v_{-x+N-m-2}u_m\xi$ later, 
we may assume $u_m\xi\not=0$ and so $Q\geq 1$. 
Let us consider a $Q\times Q$-matrix 
$$A:=((-1)^{h-j+N}\binom{-x+h+N}{j})_{h,j=0,...,Q-1} $$
consisting of coefficients of $(-1)^{x}v_{-x-2+N-j-m}(u_{m+j}\xi)$. 
It is easy to see $\det A=\pm 1$ since $\displaystyle{\binom{s+1}{j}-\binom{s}{j}=\binom{s}{j-1}}$. 
Therefore, there are polynomials $\lambda_h^m(x)\in \C[x]$ for $0\leq h< Q$ such that 
$$\begin{array}{l}
{\displaystyle (-1)^xv_{-x-2+N-m}u_m\xi}\vspace{2mm}\cr
\mbox{}\qquad{\displaystyle \equiv \sum_{h=0}^{Q-1}\lambda_h^m(x)\left( 
\sum_{j=0}^{N+h+m+2}\binom{-x+h+N}{j}(-1)^ju_{-x+h+m+N-j}(v_{-2-m-h+j}\xi)\right)}.
\end{array}$$ 
Since the coefficients of the right side at $\alpha_{-x-1}\1$ are all in $\CF_k$ by 
Lemma \ref{ID}, the above equation implies that 
a function defined by $v_{-x-2+N-m}u_m\xi$ is in $T(\CF_k)$ for  
any $v\in V^{(p-k)}$, $u\in V^{(k)}$, $\xi\in V^{(1)}$ and $m\in \Z$. 
Since $V^{(1+k)}$ is a simple $V^{(0)}$-module, $V^{(1+k)}$ is spanned by elements with 
the form $u_m\xi$ with $u\in V^{(k)}$, $\xi\in V^{(1)}$ and $m\in \Z$ and so we have  
$T(\CF_{p-k})\subseteq \CF_{k}$ for any $k$.  
Since $T^2=1$, we have the equality $T(\CF_k)=\CF_{p-k}$.  
\prend

\begin{lmm}\quad
$f(x)\in \CF_1$ is the restriction of a polynomial of $x$ on $\N$. 
\end{lmm}

\pr 
By Lemma \ref{ST}, there is $g\in \CF_0$ such that $Sg=f$. 
Since $Q_g=\{x\in \Z\mid g(x)\not=0\}$ is a finite set and $Q_g\subseteq \N$, 
we have that for $x\in \N$   
$$\begin{array}{rl}
f(x)=&{\dsp Sg(x)=(-1)^x\sum_{k=0}^x\binom{x}{k}(-1)^kg(x-k)} \cr
=&{\dsp (-1)^x\sum_{k=0}^x\binom{x}{x-k}(-1)^{x-k}g(k)=\sum_{k\in Q_g}\frac{(-1)^kg(k)}{k!}x(x-1)\cdots(x-k+1)}, 
\end{array}$$
which means that $f$ is the restriction of a polynomial on $\N$.  
\prend

\vspace{5mm}

Now we start the proof of Main Theorem. By Lemma ref{ST}, we have 
$$\begin{array}{lllllll} 
\CF_1 &  &\CF_2& &\cdots  & &\CF_{p-1} \cr
\downarrow T&\nearrow S&\downarrow T&\nearrow S &\cdots &\nearrow S&\downarrow T \cr
\CF_{p-1} & &\CF_{p-2}& & \cdots & &\CF_1, 
\end{array}$$ 
where $p$ is the order of $\sigma$. 
In particular, $T(ST)^{p-2}(\CF_1)=\CF_1$ since $S^2=T^2=1$. Therefore 
there are polynomials $f(x),g(x)\in \CF_1$ such that 
$${\dsp \sum_{j=0}^n\binom{n}{j}(p-2)^jf(n-j)=(-1)^ng(n)} \quad \mbox{ for }n\in \N.  \eqno{(4.2)}$$

\begin{lmm}\label{NF}\quad
There is no nonzero pair of polynomials satisfying (4.2). 
\end{lmm}

\pr 
This lemma was proved by Yu-ichi Tanaka and Shigeki Akiyama for rational functions 
in the case $p=3$, independently. 
The following proof is essentially given by S.~Akiyama. 

We first introduce a few notation. 
$g\in {\rm Map}(\N,\C)$ is called ``eventually positive'' if 
$\CR g(n)>0$ for a sufficiently large $n$ and 
``eventually alternating'' if $\CR g(n)\CR g(n+1)<0$ for a sufficiently large $n$, 
where $\CR g(n)$ denotes the real part of $g(n)$.  
Let $\Delta_m$ be a difference operator defined by $\Delta_m g(x):=g(x+1)-mg(x)$.
Clearly, if $g$ is eventually alternating, then so is $\Delta_m g$ for $m\geq 0$. 
We also have 
$$\Delta_{p-2}\left(\sum_{k=0}^x\binom{x}{k}(p-2)^kf(x-k)\right)=
\sum_{k=0}^x\binom{x}{k}(p-2)^kf(x+1-k) $$
by the direct calculation. 

Let us start the proof of Lemma \ref{NF}. 
Suppose false and let $(f,g)$ be a counterexample.  
We may assume that $f(x)$ is eventually positive 
(by multiplying $f(x)$ and $g(x)$ by a scalar if necessary), 
that is, there is $m\in \N$ such that 
$\CR f(k)>0$ for $k\geq m$. 
Similarly, $g(x)$ is eventually positive or negative and so 
$(-1)^xg(x)$ is eventually alternating. 
Since 
$(-1)^xg(x)=\sum_{k=0}^x\binom{x}{k}(p-2)^kf(x-k)$ is eventually alternating, so is 
$\Delta_{p-2}^m (-1)^xg(x)$. However, 
$\Delta_{p-2}^m(ST)^{p-2}f(x)=\sum_{k=0}^x\binom{x}{k}(p-2)^kf(x+m-k)$ is a 
sum of positive numbers, which is a contradiction. 
\prend

We are now able to finish the proof of Main Theorem. 
By Lemma \ref{NF}, we have $\CF_1=\{0\}$ and so $\CF_0=S(\CF_1)=\{0\}$. 
However, since $\1\in V^{(0)}, \alpha\in V^{(1)}$, 
we have $\1_{-x}\alpha=\delta_{1,x}\alpha\in \widehat{T}$ and so 
$\CF_0\not=\{0\}$, which is a contradiction. 

This completes the proof of Main Theorem. 
\prend

As last, we will prove corollaries of Main Theorem. \\

\noindent
{\bf Proof of Theorem \ref{2G}}\quad   
Let $V$ be a $C_2$-cofinite simple VOA and $G$ a finite solvable subgroup of ${\rm Aut}(V)$. 
We will prove Theorem \ref{2G} by the induction on $|G|$. 
Since $G$ is solvable, $G$ has a normal abelian subgroup $A\not=1$. 
We first assume that $G=A$ and let $1\not=\sigma\in G$ be an element of prime order. 
Then $V^\sigma$ is $C_2$-cofinite by Main Theorem. 
Furthermore, $V^\sigma$ is simple by \cite{DM}. 
Therefore, $V^G=(V^\sigma)^{G/<\sigma>}$ is also $C_2$-cofinite 
by the induction, which proves the assertion of Theorem \ref{2G}. 
So, we have $A<G$. By the minimality of $|G|$, 
$V^{A}$ is $C_2$-cofinite and it is also simple by \cite{DM}. 
Therefore, by the minimality of $|G|$, 
$V^G=(V^A)^{G/A}$ is also $C_2$-cofinite. \\

\noindent
{\bf Proof of Corollary 1}\quad   
Assume $U\cong V_L$ for some lattice $L$ and 
set $L^{\ast}=\{a\in \Q L\mid \langle a,L\rangle\subseteq \Z \}$.
We view $V$ as a $V_L$-module. 
Since $V_L$ is rational and 
the category of $V_L$-modules have a $L^{\ast}/L$-module structure, 
the actions of $G=\Hom(L^{\ast}/L,\C^{\times})$ on $V$ are induced from this structure. 
Then $V^{G}\cong U\otimes E$, 
where $E$ is a commutant of $U$ in $V$. 
By Theorem \ref{2G}, $U\otimes E$ is $C_2$-cofinite. 
If $E$ is not $C_2$-cofinite, then $E$ has a weak module $B$ containing $L(0)$-free element 
by \cite{M1} and so $U\otimes E$ has a weak module $U\otimes B$ containing $L(0)$-free elements, 
which contradicts the $C_2$-cofiniteness on $U\otimes E$. \\

\noindent
{\bf Proof of Corollary 2}\quad 
Let $U\cong L(\hf,0)$ and we view $V$ as a $U$-module. 
Since $L(\hf,0)$ is rational, $V$ is a direct sum of simple $U$-modules. 
Then $\tau$ defined by ${\rm Id}$ on $L(\hf,0)$ and $L(\hf,\hf)$ and $-{\rm Id}$ 
on $L(\hf,\st)$ as $U$-modules, 
respectively, becomes an automorphism of $V$ by the fusion rules of $L(\hf,0)$-modules.   
Then $V^{\tau}$ is $C_2$-cofinite by Main Theorem and simple by \cite{DM}.  
We then view it as a $U$-module, 
whose compositions are isomorphic to $L(\hf,0)$ or $L(\hf,\hf)$ as $U$-modules. 
Then $\sigma$ defined by ${\rm Id}$ on $L(\hf,0)$ and $-{\rm Id}$ on $L(\hf,\hf)$ as $U$-modules, 
respectively, is again an automorphism of $V^{\tau}$. 
Then $(V^{\tau})^{\sigma}=U\otimes E$, where $E$ is a commutant of $U$ in $V$. 
By Main Theorem, $U\otimes E$ is $C_2$-cofinite and so is $E$ by the same argument as above.

\end{document}